\documentclass[10pt]{amsart}
\usepackage{amssymb}
\usepackage[english,francais]{babel}
\usepackage{color}
\definecolor{darkred}{rgb}{0.8,0,0}
\definecolor{darkblue}{rgb}{0,0,0.7}
\definecolor{darkgreen}{rgb}{0,0.4,0}
\usepackage[latin1]{inputenc}
\newtheorem{theorem}{Theorem}[section]

\newtheorem{e-proposition}[theorem]{Proposition}

\newtheorem{e-definition}[theorem]{Definition\rm}

\newcommand{\RR}{\mathbb{R}}

\newcommand{\be}[1]{\begin{equation}\label{#1}}
\newcommand{\ee}{\end{equation}}

\def\og{\leavevmode\raise.3ex\hbox{$\scriptscriptstyle\langle\!\langle$~}}
\def\fg{\leavevmode\raise.3ex\hbox{~$\!\scriptscriptstyle\,\rangle\!\rangle$}}

\begin{document}
\selectlanguage{english}
\title{Selfsimilar equivalence of porous medium and $p$-Laplacian flows}

\selectlanguage{english}
\author{Ariel S\'{a}nchez}
\email{ariel.sanchez@urjc.es}
\author{Juan-Luis V\'azquez}
\email{juanluis.vazquez@uam.es}

\address{A. S.: Departamento de Matem\'atica Aplicada\\
Universidad Rey Juan Carlos\\ M\'{o}stoles. 28933 Madrid, Spain}
\address{J. L. V.:  Departamento de Matem\'
aticas\\ Universidad Aut\' onoma de Madrid\\  Cantoblanco. 28049
Madrid, Spain.}


\begin{abstract}
\selectlanguage{english}
\centerline{------------------------------------------------------------------}

We demonstrate the equivalence between the two popular models of
nonlinear diffusion, the porous medium equation and the
$p$-Laplacian equation. The equivalence is shown at the level of
selfsimilar solutions.

\vskip 0.5\baselineskip

\end{abstract}

\maketitle




\selectlanguage{english}
\section{Introduction}
\label{}

The theory of nonlinear evolution  equations of degenerate parabolic
type has been intensively studied in the last decades both for its
importance in a number of applications in mechanics and heat
propagation and also because  it combines the theory of nonlinear
evolution PDEs with geometry in the form of free boundaries. Two of
the most studied models have been the {porous medium  equation},
shortly PME, and  the $p$-Laplacian equation, shortly PLE. We write
the PME as
\begin{equation}
\label{ec.mp} u_t=\Delta(u^m/m), \quad m>1.
\end{equation}
It generalizes the heat equation, which is the linear case $m=1$,
and we also include the range of parameters $m<1$ where it is known
as the fast diffusion equation, \cite{VBook07} ($m\le 0$ is also
allowed in the general theory). On the other hand, the standard PLE
is
\begin{equation}\label{ec.plp}
u_t= \Delta_p u:= \bigtriangledown\cdot (|\bigtriangledown
u|^{p-2}\bigtriangledown u),
\end{equation}
where the parameter ranges in the interval $1<p<\infty$. The heat
equation is the case $p=2$. There has been strong interest in recent
times in the case $p=1$ as a geometrical flow, \cite{ACM}.

It was soon remarked that the theory of both equations offers
striking parallels. To mention just one, the PME has the property of
finite speed of propagation for $m>1$, while the PLE has the same
property for $p>2$. The consequence in both cases is that solutions
with compactly supported initial data $u(x,0)\ge 0$ preserve the
property of compact support w.r.t the space variable for all $t>0$,
and this implies the existence of a clear-cut interface or free
boundary separating the non-empty regions where $u>0$ and where
$u=0$. Many other  properties are known that reinforce that
parallelism, cf. \cite{dBbook}, \cite{VLN05}.

In this note we  show that the very strong connection between the
two equations can be improved into complete equivalence when we
consider special classes of solutions. Such solutions are the
self-similar solutions, i.e., solutions of one of the two  forms,
\begin{equation} \label{sim.1}
 u(x,t)=t^{-\alpha}f(x \,t^{-\beta})= t^{-\alpha}f(\eta),\qquad
 u(x,t)=(T-t)^{\alpha}f(x\,(T-t)^{\beta})=(T-t)^{\alpha}f(\eta).
\end{equation}
We will call  the solutions of the form (\ref{sim.1})-left
selfsimilar solutions of Type I, those of the form
(\ref{sim.1})-right are called Type II. The practical importance of
such solutions lies in the fact that they usually describe the
asymptotic behavior of solutions of Cauchy problem associated to
parabolic equations, as the extensive literature shows, cf.
\cite{Barbook}, \cite{VBook07} and their references. In view of the
definite progress obtained recently in the study of the PME, as
shown in the latter reference, and the key role actually played by
suitable selfsimilar solutions in the classification of the basic
qualitative and asymptotic  properties, the equivalence is of clear
significance to establish a similar theory for the PLE. We remark
that for our equations there is still a possible selfsimilarity of
Type III, with solutions of the form $u(x,t)=e^{\alpha
t}f(x\,e^{\beta t})$, but it plays a minor role in the theory. Since
the details are almost the same we will make only side mention  to
it.

A popular technique   used for  the description of these solutions
is  phase-plane analysis, cf. \cite{Barbook, MFBV, VBook07}. Our
main contribution in this Note consists in introducing a new set
of variables in the phase-plane description of the self-similar
solutions of the $p$-Laplacian equation. They are not
straightforward at all, but they have the property that the
resulting phase plane can be exactly identified with a known
version of the phase-plane for the self-similar solutions of the
PME, the only apparent difference amounting to different formulas
for the values of the parameters.
 Here is our main result.

\begin{theorem}\label{thm.kpz}
The  analysis of radial selfsimilar solutions for both the PME
(when $m\ne 1$) and the PLE (when $p\ne 2$) can be reduced to a
particular case of the autonomous ODE system
\begin{equation} \label{syst} \left\{
\begin{array}{ll}
\dot{\Psi }=\Psi \Phi\\
 \dot{\Phi}=c_1 \Phi^2-c_2 \Psi\Phi-c_3
 \Phi\pm\Psi+sgn(b)
\end{array}
\right.
\end{equation}
where the parameters $c_1$, $c_2$, $c_3$ and $b$ are explicit
functions of $n$, $m$ and $\beta$ in the PME case; of $n$, $p$ and
$\beta$ in the PLE case. The variables $\Phi$ and $\Psi$ are given
by algebraic expressions in terms of $\eta,f,f'$, which are
different for both equations. It follows that a translation rule
can be set up from the PME into the PLE and vice versa, which
involves changing the space dimension in the process. Besides,
there are in general two options for the equivalence map in either
direction.
\end{theorem}

Obviously, in the exceptional case $m=1$, $p=2$ the two equations
coincide and the identity transformation solves the equivalence
problem. We recall that there exist a large number of selfsimilar
analyses of the PME and the PLE but the equivalence has not been
remarked.

In the sequel we give whole details of these assertions: a careful
selection of the correspondence between the two sets of parameters
allows for a complete identification of the phase planes, which
implies a transformation rule for the set of selfsimilar solutions
from one equation to the other. We then derive some interesting
applications. We will also provide self-maps of the solutions of the
PME (resp. the PLE) into themselves based on an repeated use of the
transformations, using a different option at each turn. We point out
that when performing the selfsimilar analysis with radial profiles
by means of ordinary differential equations we may assume that the
space dimension is any real positive number.

\section{Phase plane analysis for the PME}

 We want to find solutions in the self-similar forms
(\ref{sim.1}) for equation (\ref{ec.mp}) posed for $x\in\RR^n$
with $m\ne 1$. This is well-known, and is described in detail in
\cite{VLN05} for instance. First, the relation between the
similarity exponents $\alpha$ and $\beta$ reads \
$(m-1)\alpha+2\beta=1$ \ for Type I selfsimilarity, and
$(m-1)\alpha+2\beta=-1$ if it belongs to Type II. Then, under the
usual assumption of radial symmetry the profile $f$ must satisfy
the ODE
\begin{equation}\label{sim.ode.mp}
\eta^{1-n}(\eta^{n-1} f^{m-1}f')'+ \alpha f + \beta \eta f'=0,
\end{equation}
where $\eta>0$.  From now on we assume that $n$ is any positive real
number.  We first introduce the variables:
\begin{equation}\label{sim.ode.mp2}
X=\eta f'/f, \qquad Y= \eta ^2 f^{1-m}.
\end{equation}
We also replace the $\eta$ variable by $r=\log \eta$. The functions
$X(r)$ and $Y(r)$ satisfy the autonomous ODE system:
\begin{equation}\label{sim.syst1.mp}
\dot{X}=(2-n)X-mX^2-(\alpha+\beta X)Y, \qquad \dot{Y}=
(2+(1-m)X)Y.
\end{equation}
where  the over-dot indicates differentiation with respect to $r$.
Now, we introduce a new pair of variables
$$
\Phi=(2+(1-m)X)/\sqrt{|b|}, \qquad \Psi=Y/|b|,
$$
where  $b=2n(m-m_c)/(m-1)$ and $m_c=(n-2)/n$. This assumes
 that $m\ne m_c$ so that $b\ne 0$. Replacing then the $r$
variable by $r_1=\sqrt{|b|}r$, so that over-dot indicates
differentiation with respect to $r_1$, the system takes the desired
quadratic form \eqref{syst}, with precise values for the constants
given by
$$
c_1=\frac{m}{m-1},\quad c_2=\beta\sqrt{|b|},\quad c_3=\frac{(n+2)(m-m_s)}{(m-1)\sqrt{|b|}},
\quad m_s=\frac{n-2}{n+2}.
$$
and the $\pm$ of the equation is $+$ for the Type I and it is $-$
for the Type II. With these values System \eqref{syst} has free
parameters $m$, $n$ and $\beta$, since $\alpha$ can be calculated
from them.  The case $m=m_c$ will be discussed below.

\section{Phase plane analysis for the PLE}

In this case, the relation between the similarity exponents
$\alpha$ and $\beta$ is \ $(p-2)\alpha+p\beta=1$ \ for Type I
selfsimilarity, and \ $(p-2)\alpha+p\beta=-1$ \ if it is of Type
II. The radially symmetric profile $f$ satisfies the  ODE
\begin{equation}\label{sim.ode}
\eta^{1-n}(\eta^{n-1} |f'|^{p-2}f')'+ \alpha  f + \beta \eta f'=0,
\end{equation}
where prime denotes differentiation with respect to $\eta>0$. From
now on we assume that $n$ is any positive real number. In a
similar way as before,  for $p\ne 2$   we introduce phase plane
variables,  a bit different from the ones in the PME case:
\begin{equation}\label{sim.ode2}
X=-\eta^2|f'|^{1-p}f', \qquad Z= \eta ^{\gamma} f,\quad
\mbox{where}\quad \gamma=\frac{p}{2-p}
\end{equation}
We thus get  the autonomous ODE system
\begin{equation}\label{sim.syst1}
\frac{p-1}{2-p}\dot{X}=-(n-\gamma)X+\alpha
Z|X|^{\frac{3-2p}{2-p}}-\beta|X|X, \qquad \ \dot{Z}= \gamma
Z-|X|^\frac{p-1}{2-p}X,
\end{equation}
where we have replaced the $\eta$ variable by $r=\log \eta$, so
that over-dot indicates differentiation with respect to $r$. This
system is not quadratic so that we perform further change with
this objective in mind. We introduce $Y=|X|^{\frac{1}{p-2}}XZ
=-\eta |f'|^{-p}f'f$ and the flow equations become
\begin{equation}\label{sim.syst2}
\dot{X}=\frac{2-p}{p-1}X\left(\gamma-n+\alpha Y-\beta |X|\right),
\qquad \ \dot{Y}= -\alpha Y^2+nY+\beta Y |X|-|X|,
\end{equation}
which is  a quadratic system if $X$ has a sign.
 For the next step we assume that $X>0$ and we set
$$
\Psi=aX, \quad \quad \Phi=-\frac{p-2}{(p-1)\sqrt{|b|}}\left(\gamma
-n +\alpha Y-\beta X\right),
$$
where $$a=\frac{1}{|b|(p-1)},\quad
b=\frac{p(n+1)(p-p_c)}{(p-2)(p-1)} \quad  p_c=\frac{2n}{n+1}.
$$
In order to proceed further we  assume that $p\ne p_c$. After all
these transformations, the flow equations become exactly the desired
(\ref{syst}). The $\pm$ of the equation is $+$ for the Type I and it
is $-$ for the Type II, and the constants are now given by
$$c_1=\frac{p-1}{p-2},\quad c_2=\beta\sqrt{|b|},\quad
c_3=\frac{(n+2)(p-p_s)}{(p-2)\sqrt{|b|}}, \quad
p_s=\frac{2n}{n+2}.
 $$
This system has free parameters $p$, $n$ and $\beta$, and $p=2$ is
excluded. We have replaced the $r$ variable by $r_1=\sqrt{|b|}r$,
so that over-dot indicates differentiation with respect to $r_1$.

\medskip

  \noindent {\bf The critical cases, $m_c$ and $p_c$.}
Some changes have to be made in the special cases  $m=m_c$ for the
PME or $p=p_c$ for the PLE, since our definitions imply that $b=0$.
we change it into $\sqrt{b}=(n-2)$ for PME and $\sqrt{b}=n$ for PLE
so that $c_3=-1$ and the independent term $sgn(b)$ disappears from
the second equation of system (\ref{syst}), which becomes
\begin{equation} \label{syst2}
\dot{\Psi }=\Psi \Phi, \qquad
 \dot{\Phi}=c_1 \Phi^2-c_2 \Psi\Phi+
 \Phi\pm\Psi.
\end{equation}
 with $c_1$ and $c_2$ as before.

\medskip

\noindent {\bf Tpe III selfsimilarity.} It offers few novelties.
The exponent relationships are $\alpha(1-m)=2\beta$ and $\alpha
(2-p)=p\beta$ respectively. We obtain a similar system, except for
the fact that the term $\pm\Psi$ disappears. Summing up, the
coefficient of $\Psi$ is $+$ for Type I, $-$ for Type II, and $0$
for Type III.

\section{Equivalence}

Identifying the two systems implies taking parameters $(m,n,\beta)$
in the first system and $(p,n',\beta')$ in the second, so that the
expressions for the free constants of System (\ref{syst}) obtain the
same values. (i) Identification of $c_1$ gives the relation $p=m+1$,
which is well-known from different dimensional considerations. (ii)
Identifying $c_3$ gives
$$
\frac{(n+2)(m-m_s)}{(m-1)\sqrt{|b|}}=\frac{(n'+2)(p-p_s)}{(p-2)\sqrt{|b'|}}
$$
Note that $b'=p(n'+1)(p-p_c)/(p-2)(p-1)$ with $p_c= 2n'/(n'+1)$.
This leads to
$$
2m(n(m-1)+2)(n'(m-1) +2(m+1))^2=(m+1)(n'(m-1)+ m+1)(n(m-1) +
2(m+1))^2,
$$
which can be written as a quadratic condition: $ A(n')^2+ Bn' +
C=0,$ with
$$
A= 2m(m-1)^2[n(m-1)+2], \quad B=-(m+1)(m-1)^3(n-2)^2, \quad
C=-(m-1)^2(m+1)^2(n-2)^2.
$$
This gives the values for $n'$ as a double function of $n$,
representing the change of dimension from PME to PLE
\begin{equation} \label{n'} n'_1=\frac{(n-2)(m+1)}{2m}, \quad
n'_2=\frac{(n-2)(m+1)}{n-2-nm}=\frac{m_c(m+1)}{m_c-m}.
\end{equation}
The double solution exists for $m\ne 0, m_c$ or $1$. The two
values of $n'$ coincide for $m=m_s$, and then $n'=n$ and $p=p_s$.
Clearly, we have
$$
\frac1{n_1'}+\frac1{n_2'}=-\frac{B}{C}=\frac{1-m}{m+1}=\frac{2-p}{p}
$$
which expresses a possible self-map of the PLE  (for fixed $p$) by
change of dimension. The situation is similar in the other sense:
 if $p\ne p_c(n')=2n'/(n'+1)$, there are
two possible values of $n$ for every given $n'$
 and we have the relationship
$$
\frac1{n_1-2}+\frac1{n_2-2}=\frac{1-m}{2m}.
$$
which expresses a possible self-map of the PME. A precedent of
self-relations can be found in King's \cite{Ki93}.

\noindent iii) Identifying $c_2$ implies that
$\beta\sqrt{|b|}=\beta'\sqrt{|b'|}$, so that
$$
\frac{\beta'^2}{\beta^2}=\frac{(m-1)p(n'+1)(p-p_c)}{2n(m-m_c)(p-2)(p-1)}
$$
which gives the corresponding two values for $\beta'$ in terms of
$\beta$
$$
\frac{\beta'_1}{\beta}=\frac{2m}{m+1}=\frac{n-2}{n'_1}, \quad
\frac{\beta'_2}{\beta}=\frac{n(m-1)+2}{m+1}=\frac{2-n}{n'_2}
$$
Note that in both cases $\beta^2 (n-2)^2=(\beta' n')^2$.

\noindent (iv) We still have to check that $sgn(b)=sgn(b')$. We have
$$
b'/b=
\frac{p(n'+1)(p-p_c)(m-1)}{2n(p-2)(p-1)(m-m_c)}=\frac{((m-1)n'+m+1)(m+1)}{2m(n(m-1)+2)}=\frac{(n')^2}{(n-2)^2}.
$$
We have used that $p=m+1$, $p_c=2n'/(n'+1)$, and
$(m-1)n'+m+1=2m(n(m-1)+2)(n')^2/(m+1)(n-2)^2 $ given by the
quadratic condition for $n'$.

\medskip

 \noindent {\bf Some general conclusions.}  (1) The
transformation is not applicable for $n=2$ since it implies $n'=0$
for both branches. Besides, the study is different for the cases
$n<2$ and
 $n>2$.

 \noindent  (2)   In the case  $n=1$ we have same
dimension $n'_2=1$ for all $m$ and then $\beta'=\beta$.  This is
well known since formally differentiating the one-dimensional PLE
gives the PME.  For the other branch corresponding to $n=1$ we have
$$
n'_1= - \frac{m+1}{2m}, \quad \beta'=-\frac{\beta}{n'_1},
$$
so that $n'_1>0$ if $-1\le m<0$, a case of very fast diffusion
studied in \cite{ERV}, see also \cite{VLN05}.  Also, $n'_1\ge 1$ if
$-1/3\le m<0$.

\noindent (3) { Cases $n\ge 3$. } (i) We have $n'=n$ if and only
if $m=m_s$ in both branches. This is the Yamabe case and it allows
us to get an explicit equivalent of the Loewner-Nirenberg solution
of the fast Diffusion Equation for the $p$-Laplacian equation.

\noindent  (ii) In the case of the first branch, $n'>0$ if $m>0$ or
$m<-1$ (a case we may disregard since then $p<0$).  In the first
branch $n'$ decreases with $m$ from $\infty$ as $m\to 0$ to
$(n-2)/2$ as $m\to\infty$.

\noindent  (iii) In the second branch $n'$ is positive for
$-1<m<m_c=(n-2)/n$ and increases from $0$ to $\infty$ in that
interval; more precisely, $n'=1$ for $m=0$

\noindent  (4) Identification of the critical cases $m=m_c$ in the
PME and $p=p_c$ in the PLE takes place for $n'=n-1$ and
$\beta'=\beta(n-2)/n'=\beta(n-2)/(n-1)$.

\noindent  (5) In the limit $m\to 1$, $p\to 2$, both equations
reduce to the linear heat equation (there is also the possibility of
converging to the eikonal equation $u_t=|\nabla u|^2$ with a
different limit process, cf. \cite{AV}). Now, extrapolating our
transformation to the linear case, $m=1$, $p=2$, we have $n'_1=n-2$
with $\beta'= \beta$, and $n_2'=2-n$ with $\beta'=\beta$.
\normalcolor

 \section{Application. Some explicit solutions}

Let us now  analyze some self-similar solutions given by the plane
of phases for the equation of porous
 medium, and  compare these solutions with the corresponding
 solution of the $p$-Laplacian equation. We will be mainly interested
 in dimensions $n\ge 3$ but lower dimensions appear in the final
 results.

\noindent {\bf Barenblatt solution and straight lines.} Consider
the equation of the trajectories with $sgn(b)=1$ and
selfsimilarity of Type I:
\begin{equation}
\label{ecua-trayec} \frac{d\Phi}{d\Psi}=\frac{c_1 \Phi^2-c_2
\Psi\Phi-c_3
 \Phi+\Psi+1}{\Psi\Phi}.
\end{equation}
We want to find the linear solutions to (\ref{ecua-trayec}).
 In order to do that we assume that  $\Phi=a_1\Psi+a_2$ and insert this form into (\ref{ecua-trayec}),  obtaining    the following
 conditions: $a_1=a_2=c_2/(c_1-1)$, and  $c_2$ must be a solution of the equation
  $c_1c_2^2-(c_1-1)c_3c_2+(c_1-1)^2=0$. Since  $c_2=\beta\sqrt{|b|}$,
  we get two values for $\beta$ and
$\beta'$ which are
\begin{eqnarray}
   \beta_1&=&\frac1{n(m-1)+2}\qquad  \mbox{  or  } \qquad \beta_2=\frac1{2m}\qquad \mbox{(porous medium)}\\
   \beta_1'&=&\frac1{n'(p-2)+p}\qquad  \mbox{  or  }  \qquad \beta_2'=\frac1{p}\qquad \mbox{
   ($p$-Laplacian)  }.
\end{eqnarray}
The profiles corresponding to the trajectories for $\beta_1$ and
$\beta_1'$ produce  the so-called Barenblatt solution, both for the
PME and the PLE. The other case is as follows: $\beta_2$ corresponds
to the dipole solution of the porous medium,
$$  f=\eta^{-\frac{n-2}{m}}\left(K-\frac{1}{b}\eta^{(mn-n+2)/m}\right)_+^{\frac1{m-1}}
$$
cf. \cite{VBook07}, Section 4.6, and $\beta_2'$ corresponds
$\alpha_2'=0$, and a profile whose derivative  is
$$ f'(\eta)=\eta^{-(n-1)/(p-1)}\left (
c-\frac{1}{(p-1)b}\eta^{(p-2)b/p}\right)_+^{1/(p-2)}
$$

\noindent {\bf Yamabe case.} It corresponds to
$p_s=m_s+1=2n/(n+2)$ (cf. \cite{VLN05}, Section 7.5 and Appendix
AIII). Choosing Type II and putting $\beta=0$, so that the
self-similar solution is a separate variables solution, we have
 $c_1=-(n-2)/4$, $c_2=c_3=0$ and $sgn(b)=1$. The  equation of the trajectories (\ref{ecua-trayec})
  turns out to be
  $$
  \frac{d\Phi}{d\Psi}=\frac{c_1\Phi^2-\Psi+1}{\Psi\Phi}.
  $$
Solving the above Bernoulli's equation,  we obtain that  the curve
$\Psi=n(1/(n-2)-1/4\Phi^2)$ is one of its solutions. The explicit
profiles that we get after  undoing all the changes of variables
are:
\begin{eqnarray}
   f(\eta)&=&(k_1+\frac{\eta^2}{4nk_1})^{-\frac{n+2}2} \qquad\mbox{  (Loewner-Nirenberg profile for the PME with $m=m_s$)  } \\
   f(\eta)&=&C(1+k_2\eta^\frac{2n}{n-2})^{-\frac{n}2}, \ C=\frac{4n}{n+2}
   \left(\frac{4n^3k_2}{n^2-4}\right)^{(n-2)/4} \quad \mbox{  (new $p$-Laplacian profile for $p=p_s$)  }
\end{eqnarray}
where $k_1, k_2>0$ are arbitrary.

\medskip

\noindent {Note.} A detailed account of the results and some
applications will appear as a separate publication \cite{ISV07}.





\end{document}